\documentclass[12pt]{article}
\usepackage{amsfonts}
\usepackage{mathrsfs}
\usepackage{amsmath,amssymb}
\baselineskip 5.5pt \lineskip 5.5pt \numberwithin{equation}{section}

\def\ad{\mbox{ad}}
\def\QED{\hfill$\Box$\par}
\def\a{\alpha}

\def\b{\beta}
\def\d{\sigma}
\def\ss{\sigma}
\def\eps{\varepsilon}

\def\DD{\mathcal {D}}

\def\g{\gamma}

\def\SS{\mathfrak{S}}
\def\HH{H}
\def\VV{\mathfrak{V}}
\def\UU{\mathcal {U}}

\def\LL{\mathscr{L}}

\def\si{\sigma}

\def\cl{\centerline}
\def\DD{{\cal D}}

\def\lll{\leftline}

\def\rar{\longrightarrow}

\def\vs{\vspace*}

\def\Inn{{\rm Inn}}
\def\Aut{{\rm Aut}}

\def\C{\mathbb{C}}

\def\Z{\mathbb{Z}}
\def\adddot{$\!\!\!${\bf.}\ \ }

\newtheorem{theo}{Theorem}[section]
\newtheorem{coro}[theo]{Corollary}
\newtheorem{lemm}[theo]{Lemma}

\begin{document}
\baselineskip 18pt

\cl{\large\bf{The derivation algebra and automorphism group}}
\cl{\large\bf{of the twisted Schr\"{o}dinger-Virasoro
algebra}\footnote {Supported by NSF grants 10471091, 10671027 of
China, ¡°One Hundred Talents Program¡± from University of Science
and Technology of China.}}\vs{6pt}

\cl{Junbo Li$^{*,\dag)}$, Yucai Su$^{\ddag)}$} \cl{\small
$^{*)}$Department of Mathematics, Shanghai Jiao Tong University,
 Shanghai 200240, China}
\cl{\small $^{\dag)}$Department of Mathematics, Changshu Institute
of Technology, Changshu 215500, China} \cl{\small
$^{\ddag)}$Department of Mathematics, University of Science and
Technology of China, Hefei 230026, China} \cl{\small E-mail:
sd\_junbo@163.com, ycsu@ustc.edu.cn} \vs{6pt}

\noindent{\bf{Abstract.}} {In this article, we determine the
derivation algebra Der$\LL$ and the automorphism group ${\rm
Aut}\LL$
of the twisted Schr\"{o}dinger-Virasoro algebra $\LL$}.
\noindent{\it{MR(2000) Subject Classification}: 17B05, 17B40, 17B65,
17B68.}\vs{28pt}\vs{6pt}

\lll{\bf1. \ Introduction}
\setcounter{section}{1}\setcounter{theo}{0}

It is well known that the infinite-dimensional Schr\"{o}dinger Lie
algebras and Virasoro algebra play important roles in many areas of
mathematics and physics (e.g., statistical physics). The original
Schr\"{o}dinger-Virasoro Lie algebra was introduced in \cite{H1}, in
the context of non-equilibrium statistical physics, containing as
subalgebras both the Lie algebra of invariance of the free
Schr\"{o}dinger equation and the Virasoro algebra. The
infinite-dimensional Lie algebra discussed in this article called
the twisted Schr\"{o}dinger-Virasoro algebra is the twisted
deformation of the original Schr\"{o}dinger-Virasoro Lie algebra.

Both original and twisted Schr\"{o}dinger-Virasoro Lie algebras are
closely related to the Schr\"{o}dinger Lie algebras and the Virasoro
Lie algebra (see \cite{H2}, \cite{HU}, \cite{ZM} and \cite{S}). They
should consequently play a role akin to that of the Virasoro Lie
algebra in two-dimensional equilibrium statistical physics.

Motivated by the research for deformations and central extensions of
both original and twisted Schr\"{o}dinger-Virasoro Lie algebras, the
sets of generators provided by the cohomology classes of the
cocycles were presented in \cite{RU}. In \cite{U}, the author
constructed vertex algebra representations of the
Schr\"{o}dinger-Virasoro Lie algebras out of a charged symplectic
boson and a free boson with its associated vertex operators. The
purpose of this article is to determine the derivation algebra and
the automorphism group of the twisted Schr\"{o}dinger-Virasoro Lie
algebras ${\LL}$. We show that the derivation algebra of ${\LL}$ is
the direct sum of three linear independent outer derivations and the
inner derivation algebra. Finally, we characterize the automorphism
group ${\rm Aut}\LL$ of $\LL$.

Now we give the definition of the Lie algebra $\LL$. A Lie algebra
$\LL$ is called a {\it twisted Schr\"{o}dinger-Virasoro Lie algebra}
(see \cite{RU}), if $\LL$ has the $\C$-basis
$$\{L_n,Y_n,M_n,C\,|\,n\in \Z\}$$
with the Lie brackets
\begin{eqnarray}
\!\!\!&\!\!\!&
[L_n,L_{n'}]=(n'-n)L_{n+n'}+\delta_{n,-n'}\frac{n^3-n}{12}C,\label{1eq2.1}
\\[4pt]\!\!\!&\!\!\!&
[L_n,Y_m]=(m-\frac{n}{2})Y_{n+m},\label{1eq2.2}\\[4pt]
\!\!\!&\!\!\!& [L_n,M_p]=pM_{n+p},\label{1eq2.3}\\[4pt]
\!\!\!&\!\!\!&[Y_m,Y_{m'}]=(m'-m)M_{m+m'},\label{1eq2.4}\\[4pt]
\!\!\!&\!\!\!&[Y_m,M_p]=[M_n,M_p]=[\LL,C]=0.\label{1eq2.55}
\end{eqnarray}
The twisted Schr\"{o}dinger-Virasoro Lie algebra has an
infinite-dimensional {\it twisted Schr\"{o}dinger subalgebra}
denoted by $\SS$ with the $\C$-basis $\{Y_n,M_n\,|\,n\in \Z\}$ and a
{\it Virasoro subalgebra} denoted by $\VV$ with the $\C$-basis
$\{L_n,C\,|\,n\in \Z\}$. The action of the Virasoro subalgebra on
the Schr\"{o}dinger subalgebra is natural. The center of $\LL$,
denoted by $C(\LL)$, is two-dimensional, spanned by $\{M_0,C\}$.
Introduce a ${\mathbb Z}$-gradation on $\LL$ by ${\rm deg}L_n={\rm
deg}Y_n={\rm deg}M_n=n$, ${\rm deg}C=0$ and decompose $\LL$ with
respect to the following gradation:
$$\LL=\bigoplus_{n\in{\mathbb Z}}\LL_{n}, \
{\LL}_0={\rm Span}_\C\{L_0,Y_0,M_0,C\},\,{\LL}_n={\rm
Span}_\C\{L_n,Y_n,M_n\},\ n\in\Z\!\setminus\!\!\{0\}.$$

Throughout the article, we denote by $\Z^*$ the set of all nonzero
integers and $\C^*$ the set of all nonzero complex numbers.

\vs{22pt}\lll{\bf2. The derivation algebra ${\rm Der}\LL$ of
$\LL$}\setcounter{section}{2}\setcounter{theo}{0}\setcounter{equation}{0}

Let $V$ be a $\LL$-module. A linear map $\varphi$ from $\LL$ to $V$
is called a {\it derivation}, if for any $x, y \in \LL$, we have
$$\varphi [x, y]=x.\varphi(y)-y.\varphi(x).$$
For $v \in V$, the map $\phi: x \rightarrow x.v$ is called an {\it
inner derivation}.

Denote by Der$(\LL,V)$ the vector space of all derivations,
Inn$(\LL,V)$ the vector space of all inner derivations. Then the
first cohomology group of $\LL$ with coefficients in $V$ is
\begin{equation}\label{1st-coh}
\HH^{1}(\LL,V)={\text{Der}_{\mathbb C}(\LL,V)}/ {\text
{Inn}_{\mathbb C}(\LL,V)}.\end{equation}
 The right-hand side is also
called the space of {\it outer derivations}.

By definition, the algebra $\LL$ is a semidirect product of
$\mathbb{Z}$-graded algebras $\LL={\VV}\ltimes{\SS},$  with
${\VV}:=\bigoplus_{n\in \mathbb{Z}}\C L_n\oplus\C C$ being the
Virasoro algebra and $\SS:=\bigoplus_{n\in \mathbb{Z}}(\C
Y_n\oplus\C M_n)$ being the twisted Schr\"{o}dinger algebra.
Clearly, $\VV$ is $\mathbb{Z}$-graded by $\VV_0=\C L_0\oplus\C C$
and ${\VV}_n=\C L_n$ for $n \in \mathbb{Z}\setminus\!\!\{0\}$ while
$\SS$ is $\mathbb{Z}$-graded by $\SS_n=\C Y_n\oplus\C M_n$ for
$n\in\Z$. By \cite{HS}, the short exact sequence
$$
\{0\}\rar\SS\rar\LL\rar\LL/\SS\rar\{0\},
$$
induces a long exact sequence
\begin{eqnarray}
\label{exact2.6-} &\{0\}\rar&
 \HH^0(\LL,\SS)\longrightarrow \HH^0(\LL,
\LL)\stackrel{f}{\longrightarrow} \HH^0(\LL,\LL/\SS)\rar \nonumber\\
&&
 \HH^1(\LL,\SS)\longrightarrow \HH^1(\LL,
\LL)\longrightarrow \HH^1(\LL,\LL/\SS)\rar\cdots
\end{eqnarray}
of $\Z$-graded vector spaces. Note that $\HH^0(\LL,\SS)=\SS^\LL=\C
M_0$, while $\HH^0(\LL,\LL)=\LL^\LL\cong \C M_0+\C C$, and
$\HH^0(\LL,\LL/\SS)=(\LL/\SS)^\LL\cong \C C$. From this we see that
the map $f$ is surjective. Thus (\ref{exact2.6-}) gives the
following exact sequence
\begin{eqnarray}
\label{exact2.6} &&\{0\}\rar
 \HH^1(\LL,\SS)\longrightarrow \HH^1(\LL,
\LL)\longrightarrow \HH^1(\LL,\LL/\SS).
\end{eqnarray}
 The right-hand side of
(\ref{exact2.6}) can be computed from the exact sequence
\begin{eqnarray}\label{exact2.7}
\{0\}\longrightarrow \HH^1({\LL}/\SS,\LL/\SS)\longrightarrow
\HH^1({\LL},\LL/{\SS}) \longrightarrow \HH^1(\SS,\LL/\SS)^\LL.
\end{eqnarray}
Note that $\HH^1({\LL}/\SS,\LL/\SS)=H^1(\VV,\VV)$ is the algebra of
outer derivations of the Virasoro algebra $\VV$, which is equal to
zero by \cite{ZM}, while $\HH^1(\SS,\LL/\SS)^\LL$ can be embeded
into \linebreak ${\rm Hom}_{\UU(\VV)}(\SS/{[\SS,\SS]},\VV)$ (where
$\UU(\VV)$ is the universal enveloping algebra of $\VV$). Therefore,
by (\ref{exact2.6}) and (\ref{exact2.7}), $\HH^1(\LL,\LL)$ will
follow from the result of the computation of $\HH^1(\LL,\SS)$ and
${\rm Hom}_{\UU(\VV)}(\SS/{[\SS,\SS]},\VV)$.
\begin{lemm}\adddot
\label{lemm2.1} $\text{\rm Der}(\LL,\SS)=\text {\rm
Der}({\LL},\SS)_{0}+{\rm {Inn}}(\LL,\SS),$ where $$\text {\rm
Der}(\LL,\SS)_{0}=\{\DD\in \text {\rm
Der}(\LL,\SS)\mid\DD(\LL_{n})\subseteq \SS_{n}, \ \text{for any} \ n
\in \mathbb{Z}\}.$$
\end{lemm}
\noindent{\it Proof.} This follows immediately from (1.2) of
\cite{F}.\QED
\begin{lemm}\adddot
\label{lemm2.2}  $\DD\in \text {\rm Der}({\LL}, \SS)_0$ if and only
if\, $\DD$ is a linear map from $\LL$ to $\SS$ satisfying the
following conditions:
\begin{eqnarray*}
&{\rm (i)}&\DD(L_n)=\big(dn+d_1\big)M_n,\\
&{\rm (ii)}&\DD (M_n)=2g_0M_n,\\
&{\rm (iii)}&\DD (Y_n)=g_0Y_n,\\
&{\rm (iv)}&\DD (C)=0,
\end{eqnarray*}
for any $n\in\Z$ and some $d,\,d_1,\,g_0\in\C$.
\end{lemm}
{\it Proof.} Let $\DD$ be a linear map form $\LL$ to $\SS$
satisfying (i)--(iv). It is clear that $\DD({\LL}_n)\subseteq
\SS_n$, and it is easy to check that $\DD$ is a derivation from
$\LL$ to $\SS$.

Let $\DD\in \text {Der}({\LL},\SS)_0$. For any $n\in\Z$, assume
\begin{eqnarray}
\label{eqln1ll} \DD(L_n)=c_nY_n+d_nM_n,\ \ \ \ c_n,\,d_n\in \C.
\end{eqnarray}
Applying $\DD$ to $[L_1, L_{n}]=(n-1)L_{n+1}$ and $[L_{-1},
L_{n+1}]=(n+2)L_{n}$, we obtain
\begin{eqnarray}
&&\!\!\!\!\!\!\!\!\!\!\!\!\!
\big(d_1-nd_{n}+(n-1)d_{n+1}\big)M_{n+1}=
\big((\frac{n}{2}-1)c_1+(n-\frac{1}{2})c_{n}-(n-1)c_{n+1}\big)Y_{n+1},\label{eqnacndn11}\\
&&\!\!\!\!\!\!\!\!\!\!\!\!\!
\big((n+2)d_{n}-d_{-1}-(n+1)d_{n+1}\big)M_{n}=\big(\frac{n+3}{2}
c_{-1}-(n+2)c_{n}+(n+\frac{3}{2})c_{n+1}\big)Y_{n}.\label{eqnacndn21}
\end{eqnarray}
Comparing the coefficients of $M_{n+1},Y_{n+1}$ in
(\ref{eqnacndn11}) and $M_{n},Y_{n}$ in (\ref{eqnacndn21}), one has
\begin{eqnarray}
&&(n-1)d_{n+1}=nd_{n}-d_1,\label{eqddcc11}\\
&&(n-1)c_{n+1}=(\frac{n}{2}-1)c_1+(n-\frac{1}{2})c_{n},\label{eqddcc21}\\
&&(n+1)d_{n+1}=(n+2)d_{n}-d_{-1},\label{eqddcc31}\\
&&(n+\frac{3}{2})c_{n+1}=(n+2)c_{n}-\frac{n+3}{2}
c_{-1}.\label{eqddcc41}
\end{eqnarray}
Subtracting (\ref{eqddcc11}) from (\ref{eqddcc31}), and
(\ref{eqddcc21}) from (\ref{eqddcc41}), we respectively obtain
\begin{eqnarray}
&&d_{n+1}-d_{n}=\frac{d_1-d_{-1}}{2},\label{eqddcc36}\\
&&c_{n+1}-c_{n}=-\frac{n}{5}(c_{-1}+c_1)+\frac{1}{5}(2c_1-3
c_{-1}).\label{eqddcc46}
\end{eqnarray}
Denoting $\frac{d_1-d_{-1}}{2}=d$ and using induction on $n$ in
(\ref{eqddcc36}), one can deduce
\begin{eqnarray}
&&d_{n}=d_1+(n-1)d,\ \ \forall\ n\in\Z.\label{eqddcc866}
\end{eqnarray}
Taking $n=-1$ in (\ref{eqddcc41}), one obtains $c_0=0$. Then taking
$n=0$ in (\ref{eqddcc46}) gives
$c_1=-c_{-1}$. From this, % According to (\ref{eqddcc469}),
(\ref{eqddcc46}) can be rewritten as $c_{n+1}-c_{n}=c_{1}.$ Thus
$c_{n}=nc_{1}$. This together with (\ref{eqddcc21}) gives
\begin{eqnarray} &&c_{n}=0,\ \ \ \forall\
n\in\Z.\label{eqddcc929}
\end{eqnarray}
According to the equations (\ref{eqddcc866}) and (\ref{eqddcc929}),
for any $n\in\Z$, we can rewrite (\ref{eqln1ll}) as
\begin{eqnarray}
\label{eqln221} \DD (L_n)=\big(d_1+(n-1)d\big)M_n,\ \ {\rm{where}}\
\ d,\,d_1\in\C.
\end{eqnarray}
Re-denoting  $d_1-d$ by $d_1$, we obtain (i).

Applying $\DD$ to $[L_{-3},L_{3}]=6L_0-2C$, one has
\begin{eqnarray*}6(d_1-d)M_0-2\DD(C)\!\!&=&\!\!6\DD(L_0)-2\DD(C)=
\DD\big([L_{-3},L_{3}]\big)\\
\!\!&=&\!\![(d_1-4d)M_{-3},L_3]+[L_{-3},(d_1+2d)M_3]\\
\!\!&=&\!\!3(d_1-4d)M_{0}+3(d_1+2d)M_0\\
\!\!&=&\!\!6(d_1-d)M_0.
\end{eqnarray*}
This gives $\DD(C)=0.$ Thus (iv) follows.

For any $n\in\Z$, write
\begin{eqnarray} \label{eqln1ll1}
\DD(Y_n)=g_nY_n+h_nM_n,\ \ \ \ g_n,\,h_n\in \C.
\end{eqnarray}
Applying $\DD$ to $[L_1,Y_n]=(n-\frac{1}{2})Y_{n+1}$ and comparing
the coefficients of $Y_{n+1}$ and $M_{n+1}$, we obtain
\begin{eqnarray} \label{eqln2www2}
&&(n-\frac{1}{2})g_n=(n-\frac{1}{2})g_{n+1},\ \ \
nh_n=(n-\frac{1}{2})h_{n+1}.
\end{eqnarray}
Applying $\DD$ to $[L_{-1},Y_1]=\frac{3}{2}Y_{0}$ and comparing the
coefficients of $Y_{0}$ and $M_{0}$, we obtain $g_1=g_0,\,
h_1=\frac{3}{2}h_0.$ Then using induction on $n$ in
(\ref{eqln2www2}), one obtains
\begin{eqnarray} \label{eqldsa2}
&&g_n=g_{0},\ \ \ h_n=0,\ \ \forall\ n\in\Z.
\end{eqnarray}
Therefore, for any $n\in\Z$, one can rewrite (\ref{eqln1ll1}) as
$\DD(Y_n)=g_0Y_n,$ i.e., (iii) follows.

For any $m,n\in\Z$, applying $\DD$ to $[Y_m,Y_n]=(n-m)M_{n}$ and
using (iii), one can easily obtain (ii). \QED\vskip6pt

Set $d_1=1,\,d=g_0=0;\ \, d=1,\,d_1=g_0=0;\ \,g_0=1,\,d_1=d=0$ in
$\DD\in \text {Der}({\LL}, \SS)_0$ respectively, we obtain three
 derivations from $\LL$ to $\SS$:
\begin{eqnarray*}
\begin{array}{lllllll}
&&\DD_1:\ \DD_1(L_n)=M_n,&\DD_1(M_n)=\DD_1(Y_n)=\DD_1(C)=0;\\[7pt]
&&\DD_2:\ \DD_2(L_n)=nM_n,&\DD_2(M_n)=\DD_2(Y_n)=\DD_2(C)=0;\\[7pt]
&&\DD_3:\ \DD_3(Y_n)=Y_n,& \DD_3(M_n)=2M_n, \ \DD_3(L_n)=\DD_3(C)=0;
\end{array}\end{eqnarray*} where $n\in \mathbb{Z}$. It is straightforward to
verify that $\DD_i,\,i=1,2,3,$ are linear independent outer
derivations. Thus (\ref{1st-coh}) and Lemma \ref{lemm2.2} prove
Lemma \ref{lemm2.3} below.
\medskip
\begin{lemm}\adddot
\label{lemm2.3}
$H^1(\LL,{\SS})=\mathbb{C}\DD_1\dot{+}\mathbb{C}\DD_2
\dot{+}\mathbb{C}\DD_3$.
\end{lemm}

\begin{lemm}\adddot
\label{lemm2.4} ${\rm Hom}_{\UU(\VV)}(\SS/{[\SS,\SS]},\VV)=0$.
\end{lemm}
{\it Proof.} By the definition of $\SS$, we have $[\SS,\SS]={\rm
Span}_\C\{M_n\,|\,n\in\Z\}$. For any $f \in \text
{Hom}_{\UU(\VV)}(\SS/{[\SS,\SS]},\VV)$ and any $n\in\Z$, one can
assume
$$f(Y_n)=\mbox{$\sum\limits_{k}$}p_k^{(n)}L_{k}+c^{(n)}C,\ \ \
\mbox{for some \ }p_k^{(n)},\,\,c^{(n)}\in \mathbb{C}.$$ For any
$m,n\in\Z$, applying $f$ to both sides of $[L_{m},
Y_n]=(n-m)Y_{m+n}$, using $f([L_m,Y_n])=[L_m,f(Y_n)]$, one has
\begin{eqnarray}
\label{eqwa28}
\mbox{$\sum\limits_k$}(k-m)p_k^{(n)}L_{k+m}+p_{-m}^{(n)}\frac{m^3-m}{12}C
=(n-m)\big(\mbox{$\sum\limits_{k}$}p_k^{(m+n)}L_{k}+c^{(m+n)}C\big).
\end{eqnarray}
Comparing the coefficients of $C$ and $L_k$ in (\ref{eqwa28}), one
has
\begin{eqnarray}
&&p_{-m}^{(n)}\frac{m^3-m}{12}=(n-m)c^{(m+n)},\label{eqwa228}\\
&&(n-m)p_k^{(m+n)}=(k-2m)p_{k-m}^{(n)}.\label{eqwa268}
\end{eqnarray}
Taking $m=0,\,n\neq0$ and $m=1,n=-1$ in (\ref{eqwa228})
respectively, one immediately obtains $C^{(n)}=0$ for $n\in\Z.$
Using this in (\ref{eqwa228}), one has
\begin{eqnarray}
&&p_{k}^{(n)}=0\ \ \ \mbox{\rm for all\ \,}\ k\notin
\{-1,0,1\}.\label{eqwa2928}
\end{eqnarray}
If $k\in\{-1,0,1\}$, by choosing some $m$ such that $2n-m\ne0$ and
$k-m\notin\{-1,0,1\}$ and replacing $n$ by $n-m$ in (\ref{eqwa268}),
and using (\ref{eqwa2928}), one again obtains $p_k^{(n)}=0$. This
proves $f(Y_n)=0$ for all $n\in\Z.$\QED

\medskip
\begin{theo}\adddot
\label{theo2.5} $\text {Der}{\LL}=\mathbb{C}\DD_1\dot{+}
\mathbb{C}\DD_2\dot{+}\mathbb{C}\DD_3\dot{+}\text {ad}\LL$.
\end{theo}
{\it Proof.} By Lemma \ref{lemm2.4}, $\HH^1(\SS,\LL/\SS)^{\LL}=0$.
According to the exact sequence (\ref{exact2.7}) and the fact
$\HH^1(\LL/\SS,\LL/\SS)=0$, we have $\HH^1(\LL,\LL/{\SS})=0$. So by
Lemma \ref{lemm2.3} and
 the exact sequence (\ref{exact2.6}), one has
 $$\HH^1(\LL,\LL)=
 \mathbb{C}\DD_1\dot{+}\mathbb{C}\DD_2\dot{+}\mathbb{C}\DD_3.$$
 Hence the theorem follows from $\HH^1(\LL,\LL)
 =\text {Der}\LL/\text {ad}\LL.$\QED

\vs{18pt} \lll{\bf3. The automorphism group ${\rm Aut}\LL$ of
$\LL$}\setcounter{section}{3}\setcounter{theo}{0}\setcounter{equation}{0}

Throughout this section, we denote by $\Inn(\LL)$ the set of inner
automorphisms. Then $\Inn(\LL)$ is a normal subgroup of $\text
{Aut}{\LL}$ and \begin{equation} \label{inn-group} \mbox{
$\Inn(\LL)$ is generated by $\text {exp}(k\text{ad} L_{0}),\,\,\text
{exp}(k\text{ad} Y_{n})$ and $\text {exp}(k\text{ad}
M_{n})$},\end{equation}
 where $n \in
\mathbb{Z},k\in\C$.

Denote $J_1=\text {Span}_{\mathbb{C}}\{M_n,C\,|\, n\in \mathbb{Z}\}$
and $J_2=\text {Span}_{\mathbb{C}}\{Y_n,M_n,C\,|\,n\in
 \mathbb{Z}\}$. They both are ideals of $\LL$. And the following lemma
 is clear.
\begin{lemm}\adddot
\label {Lemma 3.1} For any $\ss\in \text {Aut}{\LL}$, one has
$\ss(J_i)\subseteq J_i\,(i=1,2)$ and $\ss\big(C(\LL)\big)\subseteq
C(\LL)$.
\end{lemm}

For any $a\in\C^*$, we can define the following maps on $\LL$.
\begin{eqnarray}
&\ss_a&\!\!\!: \ \ \ L_n\rar a^nL_n,\ \ \ Y_n\rar a^nY_n,\ \ \
M_n\rar a^nM_n,\ \ \
C\rar C;\label{eqssa}\\
&\varepsilon&\!\!\!: \ \ \ L_n\rar -L_{-n},\ \ Y_n\rar -Y_{-n},\ \
M_n\rar -M_{-n},\ \ C\rar -C;\label{eqvar}\\
&\psi_a&\!\!\!: \ \ \ L_n\rar L_n,\ \ \ \ \ \ Y_n\rar aY_n,\ \ \ \
\, M_n\rar a^2M_n,\ \ \ C\rar C.\label{eqssa'}
\end{eqnarray}
It is easy to see that
$\ss_a\in{\rm{Inn}}\LL,\,\,\eps,\psi_a\in\rm{Aut}\LL\setminus\rm{Inn}\LL$.
\begin{lemm}\adddot
\label{Lemma 3.2} For any $\ss\in \Aut({\LL})$, there exists an
automorphism $\xi\in\Inn({\LL})$ of the form ${\rm
exp}\big(\sum_{j\neq0}b'_{j}{\rm ad}Y_{j}+\sum_{k\neq0}c'_{k}{\rm
ad}M_{k}\big)$ such that
$\xi^{-1}\d(L_0)=a(L_0+b_0Y_0+c_0M_0+d_0C)$, where $a\in
\C^*,b_0,c_0,d_0\in \mathbb{C}$.
\end{lemm}
{\it Proof.} For any $\ss\in {\rm{Aut}}(\LL)$, we can assume that
\begin{eqnarray}
\label{eqd011} &&\d(L_0)=\mbox{$\sum\limits_{i\neq0}$}a_{i}L_{i}
+\mbox{$\sum\limits_{j\neq0}$}b_{j}Y_{j}+\mbox{$\sum\limits_{k\neq0}$}c_{k}M_{k}+
a(L_0+b_0Y_0+c_0M_0+d_0C),
\end{eqnarray}
where $a_{i},b_{j},c_{k}\in \mathbb{C},\,a\in\C^*$ (we must have
$a\ne0$, otherwise, $\si(L_0)$ would be $\ad$-locally nilpotent
while $L_0$ is $\ad$-semi-simple).

We shall find some $\xi={\rm exp}\big(\sum\limits_{j\neq0}b'_{j}{\rm
ad}Y_{j}+\sum\limits_{k\neq0}c'_{k}{\rm ad}M_{k}\big)\in {\rm
Inn}(\LL)$, such that
\begin{eqnarray}
\label{d L_0=rm exp22} &&\d(L_0)={\rm
exp}\big(\mbox{$\sum\limits_{j\neq0}$}b'_{j}{\rm
ad}Y_{j}+\mbox{$\sum\limits_{k\neq0}c'_{k}$}{\rm
ad}M_{k}\big)\big(a'(L_0+b_0'Y_0+c_0'M_0+d_0'C)\big).
\end{eqnarray}
for some $a'\in \C^*,b_j',c_k',d_0'\in \mathbb{C}$. We only need to
solve $a',b'_j,c'_k,d'_0$ by comparing the coefficients of
$L_{i},Y_{j},M_{k}$ and $C$, which gives
\begin{eqnarray}
&&a=a',\ b_0=b_0',\ c_0=c_0',\ d_0=d_0',\ b_{j}'=-a^{-1}j^{-1}b_{j},\label{equanssw11}\\
&&c'_{k} =-b_0b'_{k}-a^{-1}k^{-1}c_{k}
-\frac{\sum\limits_{i+l=k}l(l-i)k^{-1}a^{-2}i^{-1}l^{-1}b_{i}b_{l}}{2}.\label{equanssw22}
\end{eqnarray}
Applying $\xi^{-1}$ to both sides of (\ref{eqd011}) and using
(\ref{equanssw11}), we obtain
\begin{eqnarray}
\label{d L_0=rm exp12}
&&\xi^{-1}\d(L_0)=a(L_0+b_0Y_0+c_0M_0+d_0C).
\end{eqnarray}
\begin{lemm}\adddot
\label{Claim1}
$b_0=0$.
\end{lemm}
{\it Proof.} Since $L_0$ is an $\ad$-semi-simple element in $\LL$,
then $\xi^{-1}\d(L_0)=a(L_0+b_0Y_0+c_0M_0+d_0C)$ is also
$\ad$-semi-simple in $\LL$. By linear algebra, the following matrix
(\,given by the adjoint action of
${\rm{ad\,}}a(L_0+b_0Y_0+c_0M_0+d_0C)$ on the basis
$\{L_n,Y_n,M_n,C\}$ of $\LL_n$\,)
$${\left(\begin{array}{cccc}
an\ \ \ \ \ {0} \ \ \ \ \ {0}\ \ \ \ {0}\\
\!\!{\frac{ab_0n}{2}}\ \ \ \,{an} \ \ \ \,\,{0}\ \ \ \,\,{0}\\
\ \,{0}\ \ \,\,\ {ab_0n} \ \,\,{an}\ \ {0}\\ \ \ \,\!\!{0} \,\ \ \ \
\ {0} \ \ \ \ \ {0}\ \ \ \ {0}
\end{array}\right)}$$
can be diagonalized, which forces $b_0=0$. This lemma follows.\QED

According to Lemma \ref{Lemma 3.2} and \ref{Claim1}, for any
$\ss\in{\rm{Aut}}(\LL)$, by replacing $\ss$ with $\xi\ss$ for some
$\xi\in{\rm{Inn}}(\LL)$, we can write
\begin{eqnarray}\label{eqssL02}
\ss(L_0)=a(L_0+c_0M_0+d_0C),
\end{eqnarray}
where $a\in \C^*,c_0,d_0\in \mathbb{C}$.
\begin{lemm}\adddot
\label{Lemma 3.3} For any $\ss\in{\rm{Aut}}(\LL)$, by replacing
$\ss$ with $\eps\ss$ if necessary where $\eps$ is defined in
$(\ref{eqvar})$, we can write
\begin{eqnarray*}
\ss(L_0)=L_0+c_0M_0+d_0C,\,
\ss(Y_n)=e'_nY_{n}+e''_{n}M_{n}+\delta_{n,0}e_0C,\,\ss
(M_n)=f'_{n}M_{n}+\delta_{n,0}f_0C,
\end{eqnarray*}
for some $n\in \mathbb{Z},\,e'_{n},\,f_n'\in
\mathbb{C}^*,\,e_0,f_0,e_n''\in \mathbb{C}$.
\end{lemm}
{\it Proof.} By Lemma \ref{Lemma 3.1}, we can write
\begin{eqnarray}
&&\ss(Y_n)=\mbox{$\sum\limits_{i}$}e'_iY_{i}+\mbox{$\sum\limits_{j}$}e''_{j}M_{j}+e_nC,\label{eqdym1}\\
&&\ss(M_n)=\mbox{$\sum\limits_{i}$}f'_{i}M_{i}+f_nC,\label{eqdmp1}
\end{eqnarray}
where
$e'_i\in\C^*,f'_i\in\C^*\,(i\neq0),\,e''_j,e_n,f'_0,f_n\in\C,\,n\in
\mathbb{Z}$.

Applying $\ss$ to $[L_0,Y_n]=nY_n$ and $[L_0,M_n]=nM_n$, we obtain
\begin{eqnarray*}
&&\Big[a(L_0+c_0M_0+d_0C),
\mbox{$\sum\limits_{i}$}e'_iY_{i}+\mbox{$\sum\limits_{j}$}e''_{j}M_{j}+e_nC\Big]
=n(\mbox{$\sum\limits_{i}$}e'_iY_{i}+\mbox{$\sum\limits_{j}$}e''_{j}M_{j}+e_nC),\\
&&\Big[a(L_0+c_0M_0+d_0C),
\mbox{$\sum\limits_{i}$}f'_{i}M_{i}+f_nC\Big]=n(\mbox{$\sum\limits_{i}$}f'_{i}M_{i}+f_nC).
\end{eqnarray*}
That is,
\begin{eqnarray}
&&\mbox{$\sum\limits_{i}$}(ai-n)e'_iY_{i}+\mbox{$\sum\limits_{j}$}(aj-n)e''_{j}M_{j}
-ne_nC=0,\label{eq3.38.1}\\
&&\mbox{$\sum\limits_{i}$}(ai-n)f'_{i}M_{i}-nf_nC=0.\label{eq3.39.1}
\end{eqnarray}
Comparing the coefficients of $Y_{i}$, $M_{j}$ and $C$ in the above
two equations, we obtain
\begin{eqnarray}\label{iajen2}
&&i=\frac{n}{a}\in\Z,\,\ (aj-n)e''_{j}=0,\ \ e_n=\delta_{n,0}e_0\ \
{\rm and}\ \ \,f_n=\delta_{n,0}f_0,
\end{eqnarray}
for all $n\in\Z$. It follows that
\begin{eqnarray}\label{ai=pm1}
a=\pm1,\ \ i=an \ \ \ {\rm{and}}\ \ \ e_j''=0\ \ \ {\rm{if}}\ \
j\neq an.
\end{eqnarray}
According to (\ref{iajen2}) and (\ref{ai=pm1}), one can respectively
rewrite (\ref{eqdym1}) and (\ref{eqdmp1}) as
\begin{eqnarray}
&&\ss(Y_n)=e'_{an}Y_{an}+e''_{an}M_{an}+\delta_{n,0}e_0C,\label{eqdym12}\\
&&\ss(M_n)=f'_{an}M_{an}+\delta_{n,0}f_0C.\label{eqdmp12}
\end{eqnarray}
For one case $a=1$, the lemma is right. For the other case $a=-1$,
replacing $\ss$ with $\eps\ss$, we can rewrite (\ref{eqssL02}),
(\ref{eqdym12}) and (\ref{eqdmp12}) as those in the lemma. Then this
lemma follows.\QED
\begin{lemm}\adddot
\label{Lemma 3.4} Let $\ss\in \Aut({\LL})$ be such that
\begin{eqnarray*}
\ss(L_0)=L_0+c_0M_0+d_0C,\,
\ss(Y_n)=e'_nY_{n}+e''_{n}M_{n}+\delta_{n,0}e_0C,\,\ss
(M_n)=f'_{n}M_{n}+\delta_{n,0}f_0C,
\end{eqnarray*}
where $n\in\mathbb{Z},\,e'_{n},\,f_n'\in
\mathbb{C}^*,\,e_0,f_0,e_n''\in \mathbb{C}$. Then by replacing $\ss$
by $\psi_x\ss_y\ss$ for some $x,y\in\Z^*$ $($where $\ss_y$ and
$\psi_x$ are defined in $(\ref{eqssa})$ and $(\ref{eqssa'})$
respectively$)$, one can suppose $\ss=\delta_{\a,\b,\g}$ for some
$\a,\b,\g\in\C$, where $\delta_{\a,\b,\g}$ is defined by
\begin{eqnarray}
\delta_{\a,\b,\g}
(C)\!\,&=&\!C,\label{eqLateq4200}\\
\delta_{\a,\b,\g}
(M_n)\!\!\!&=&\!M_{n},\label{eqLateq3200}\\
\delta_{\a,\b,\g} (Y_n)\!&=&\!Y_{n}+2\a nM_{n},\label{eqLateq2200}\\
\delta_{\a,\b,\g}( L_n)\!\!&=&\!L_{n}+\a nY_{n}+\big(\a^2n^2+\b
n+\g\big)M_{n}.\label{eqLateq1200}
\end{eqnarray}
\end{lemm}
{\it Proof.} For any $n\in\Z^*$, we can write
\begin{eqnarray}
\label{eqdln11} &&\d(L_n)=\mbox{$\sum\limits_{i}$}u_{i}L_{i}
+\mbox{$\sum\limits_{j}$}v_{j}Y_{j}+\mbox{$\sum\limits_{k\neq0}$}w_{k}M_{k}+C_n,
\end{eqnarray}
where $u_i\in\C^*,v_j,w_k\in \mathbb{C},\ C_n\in C(\LL)$.

For any $n\in\Z^*$, applying $\d$ to $nL_{n}=[L_{0},L_{n}]$, we
obtain
\begin{eqnarray*}
\mbox{$\sum\limits_{i}$}nu_{i}L_{i}
+\mbox{$\sum\limits_{j}$}nv_{j}Y_{j}+\mbox{$\sum\limits_{k\neq0}$}nw_{k}M_{k}+nC_n
\!\!\!&=&\!\!\! n\d(L_{n})=
\d\big([L_{0},L_{n}]\big)\\
\!\!\!&=&\!\!\!\mbox{$\sum\limits_{i}$}iu_{i}L_{i}
+\mbox{$\sum\limits_{j}$}jv_{j}Y_{j}+\mbox{$\sum\limits_{k\neq0}$}kw_{k}M_{k}.
\end{eqnarray*}
Comparing the coefficients of $L_i,\,Y_i,\,M_i$ and $C_n$, one has
\begin{eqnarray*}
u_i=v_i=w_i=0 \ \,{\rm{if}}\ \,i\neq n,\ \,{\rm{and}}\ \,C_n=0\ \
\,{\rm{for\ \,any}}\ \,n\in\Z^*.
\end{eqnarray*}
Hence (\ref{eqdln11}) can be rewritten as
\begin{eqnarray}
\label{eqdln22}
&&\d(L_n)=u_{n}L_{n}+v_{n}Y_{n}+w_{n}M_{n} \ \
{\rm{where}}\ \ n\in\Z^*.
\end{eqnarray}
For $n\neq0,\pm1$, applying $\d$ to $(n-1)L_{n+1}=[L_{1},L_{n}]$, we
obtain
\begin{eqnarray*}
(n-1)\si(L_{n+1})\!\!\!&=&\!\!\!(n-1)(u_{n+1}L_{n+1}+v_{n+1}Y_{n+1}+w_{n+1}M_{n+1})=[\si(L_1),\si(L_n)]\\
\!\!\!&=&\!\!\!(n-1)u_{1}u_{n}L_{n+1}+\big((n-\frac{1}{2})u_{1}v_{n}+(\frac{n}{2}-1)v_{1}u_{n}\big)Y_{n+1}\\
&&+\big(nu_{1}w_{n}+(n-1)v_{1}v_{n}-w_{1}u_{n}\big)M_{n+1}.
\end{eqnarray*}
Comparing the coefficients of $L_{n+1},\,Y_{n+1}$ and $M_{n+1}$, one
has
\begin{eqnarray}
&&u_{n+1}=u_{1}u_{n},\label{equ1}\\
&&(n-1)v_{n+1}=(n-\frac{1}{2})u_{1}v_{n}+(\frac{n}{2}-1)v_{1}u_{n},\label{eqv1}\\
&&(n-1)w_{n+1}=nu_{1}w_{n}+(n-1)v_{1}v_{n}-w_{1}u_{n}.\label{eqwuv}
\end{eqnarray}
Applying $\d$ to $[L_{1},L_{-1}]=-2L_{0}$, we obtain
\begin{eqnarray*}
-2\si(L_0)\!\!\!&=&\!\!\!-2(L_0+c_0M_0+d_0C)=\ss\big([L_{1},L_{-1}]\big)\\
\!\!\!&=&\!\!\!-2u_{1}u_{-1}L_{0}-\frac{3}{2}(u_{1}v_{-1}+v_{1}u_{-1})Y_{0}-(u_{1}w_{-1}+2v_{1}v_{-1}+w_{1}u_{-1})M_{0}.
\end{eqnarray*}
Comparing the coefficients of $L_{0},\,Y_{0}$, $M_{0}$ and $C$, one
has
\begin{eqnarray}
&&u_{-1}=u_{1}^{-1},\label{equ-11}\\
&&v_{-1}=-u_1^{-2}v_{1},\label{eqv-11}\\
&&c_0=\frac{1}{2}(u_{1}w_{-1}+2v_1v_{-1}+w_{1}u_{1}^{-1}),\ \
d_0=0.\label{eqnau-1v-1}
\end{eqnarray}
By (\ref{eqssL02}), we can rewrite $\ss L_0$ as
\begin{eqnarray*}
\ss(L_0)=L_0+\frac{1}{2}(u_{1}w_{-1}+2v_1v_{-1}+w_{1}u_{1}^{-1})M_0.
\end{eqnarray*}
Applying $\d$ to $[L_{2},L_{-2}]=-4L_{0}+\frac{1}{2}C$, we obtain
\begin{eqnarray*}
-4(L_0+c_0M_0)+\frac{1}{2}\ss C
\!\!\!&=&\!\!\!\ss\big([L_{2},L_{-2}]\big)\\
\!\!\!&=&\!\!\!-4u_{2}u_{-2}L_{0}-3(u_{2}v_{-2}+v_{2}u_{-2})Y_{0}+\frac{1}{2}u_{2}u_{-2}C-2u_{2}w_{-2}M_{0}\\
&&-4v_{2}v_{-2}M_{0}-2w_{2}u_{-2}M_{0}.
\end{eqnarray*}
Comparing the coefficients of $L_{0},\,Y_{0}$, using
(\ref{eqnau-1v-1}) and Lemma \ref {Lemma 3.1}, one has (\,noting
that $\si C\in C(\LL)$\,)
\begin{eqnarray}
&&u_{2}=u_{-2}^{-1}=u_{1}^{2},\label{equ21}\\
&&v_{2}=-u_{2}v_{-2}u_{-2}^{-1}=2u_{1}v_1,\label{eqv21}\\
&&\ss
C=C+4(u_{1}w_{-1}+2v_1v_{-1}+w_{1}u_{1}^{-1}-u_{2}w_{-2}-2v_{2}v_{-2}-w_{2}u_{-2})M_{0}.\label{equdelaC1}
\end{eqnarray}
According to (\ref{equ-11}), (\ref{equ21}) and (\ref{equ1}), one
obtains
\begin{eqnarray}
&&u_{n}=u_{1}^{n} \ \ \ {\rm{for\ \,any}}\ \,n\in\Z^*.\label{equnn1}
\end{eqnarray}
Similarly, by (\ref{eqv-11}), (\ref{eqv21}) and (\ref{eqv1}), one
has
\begin{eqnarray}
&&v_{n}=nu_{1}^{n-1}v_1 \ \ \ {\rm{for\ \,any}}\
\,n\in\Z^*.\label{eqvnn1}
\end{eqnarray}
Noting the equations (\ref{equnn1}) and (\ref{eqvnn1}), for
$n\neq0,\pm1$, we can rewrite (\ref{eqwuv}) as
\begin{eqnarray}
&&(n-1)w_{n+1}=nu_{1}w_{n}+(n-1)u_{1}^{n-1}v_1^2-w_{1}u_{1}^{n}.\label{eqwuv22}
\end{eqnarray}
For any $n\in\Z$, applying $\ss$ to both sides of
$[L_1,Y_n]=(n-\frac{1}{2})Y_{n+1}$, one has
\begin{eqnarray*}
&&(n-\frac{1}{2})(u_{1}e'_n-e'_{n+1})Y_{n+1}+\big((n-1)v_{1}e'_n+nu_{1}e''_{n}-(n-\frac{1}{2})e''_{n+1}\big)M_{n+1}\\
&&=(n-\frac{1}{2})\delta_{n+1,0}e_0C.
\end{eqnarray*}
Comparing the coefficients of $Y_{n+1},M_{n+1}$ and $C$ in the above
equation, one has
\begin{eqnarray}
&&e_0=0,\label{eqe011}\\
&&e'_{n+1}=u_{1}e'_n,\label{eqe'n+11}\\
&&(n-\frac{1}{2})e''_{n+1}=(n-1)v_{1}e'_n+nu_{1}e''_{n}.\label{eqe''n+11}
\end{eqnarray}
Using induction on $n$ in (\ref{eqe'n+11}), one has
\begin{eqnarray}
&&e'_{n}=u_{1}^{n}e'_0.\label{eqe''n2001}
\end{eqnarray}
Applying $\ss$ to both sides of $[L_{-1},Y_1]=\frac{3}{2}Y_{0}$, one
has
\begin{eqnarray*}
\frac{3}{2}\ss(Y_{0})&=&\frac{3}{2}e'_1u_{-1}Y_{0}+(2e'_1v_{-1}+e''_{1}u_{-1})M_{0}\\
&=&\frac{3}{2}e'_1u_{1}^{-1}Y_{0}+(-2e'_1u_1^{-2}v_{1}+e''_{1}u_{1}^{-1})M_{0}\\
&=&\frac{3}{2}e'_0Y_{0}+(-2u_{1}e'_0u_1^{-2}v_{1}+2v_1e'_0u_{1}^{-1})M_{0}\\
&=&\frac{3}{2}e'_0Y_{0}.
\end{eqnarray*}
So we can rewrite $\ss Y_0$ as $\ss Y_{0}=e'_0Y_{0},$ and $e''_0=0.$
Using induction on $n$ in (\ref{eqe''n+11}), one has
\begin{eqnarray}
&&e''_{n}=2nu_{1}^{n-1}v_1e'_0.\label{eqe''n21}
\end{eqnarray}
By the equations (\ref{eqe011}), (\ref{eqe''n2001}) and
(\ref{eqe''n21}), for any $n\in\Z$, we can rewrite $\ss(Y_n)$ as
\begin{eqnarray}
\ss(Y_n)=u_{1}^{n}e'_0Y_{n}+2nu_{1}^{n-1}v_1e'_0M_{n}.
\end{eqnarray}
For any $n\in\Z$, applying $\ss$ to both sides of
$[L_1,M_n]=nM_{n+1}$, one has
\begin{eqnarray}\label{eqnuf'1}
&&nu_{1}f'_{n}M_{n+1}=nf'_{n+1}M_{n+1}+n\delta_{n+1,0}f_0C.
\end{eqnarray}
For $n=-1$, comparing the coefficients of $M_{0}$ and $C$ in
(\ref{eqnuf'1}), one has
\begin{eqnarray}\label{eqnuf'22}
&&f'_{-1}=u_{1}^{-1}f'_{0},\ \ \ f_0=0.
\end{eqnarray}
For $n\neq0,-1$, using (\ref{eqnuf'22}) and comparing the
coefficients of $M_{n+1}$ and $C$ in (\ref{eqnuf'1}), one has
\begin{eqnarray}\label{eqnuf'66}
&&f'_{n+1}=u_{1}f'_{n}.
\end{eqnarray}
%For $n=-1$, comparing the coefficients of $M_{0}$ and $C$ in
%(\ref{eqnuf'1}), one has
%\begin{eqnarray}\label{eqnuf'22}
%&&f'_{-1}=u_{1}^{-1}f'_{0},\ \ \ f_0=0.
%\end{eqnarray}
%For $n\neq0,-1$, using (\ref{eqnuf'22}) and comparing the
%coefficients of $M_{n+1}$ and $C$ in (\ref{eqnuf'1}), one has
%\begin{eqnarray}\label{eqnuf'66}
%&&f'_{n+1}=u_{1}f'_{n}.
%\end{eqnarray}
Applying $\ss$ to both sides of $[L_2,M_{-1}]=-M_{1}$ and comparing
the coefficients of $M_1$, one has
\begin{eqnarray}\label{eqnf'1'1}
f'_{1}=u_{2}f'_{-1}=u_1^{2}u_{1}^{-1}f'_{0}=u_1f'_{0}.
\end{eqnarray}
According to the equations (\ref{eqnuf'22}), (\ref{eqnf'1'1}) and
using induction on $n$ in (\ref{eqnuf'66}), one has
\begin{eqnarray}\label{eqnuf'88}
&&f'_{n}=u_{1}^nf'_{0},\ \ \forall\, n\in\Z.
\end{eqnarray}
By the equations (\ref{eqnuf'22}) and (\ref{eqnuf'88}), for any
$n\in\Z$, we can rewrite $\ss(M_n)$ as
\begin{eqnarray}\label{sseqMn1}
\ss(M_n)=u_{1}^nf'_{0}M_{n}.
\end{eqnarray}
Applying $\ss$ to $[Y_{-1},Y_{1}]=2M_{0}$ and comparing the
coefficients of $M_0$, one has
\begin{eqnarray}\label{eqYYM068}
&&f'_{0}={e'_0}^2.
\end{eqnarray}
By the equations (\ref{sseqMn1}) and (\ref{eqYYM068}), for any
$n\in\Z$, we can rewrite $\ss(M_n)$ as
\begin{eqnarray}\label{sseqMn2}
\ss(M_n)=u_{1}^n{e'_0}^2M_{n}.
\end{eqnarray}
By now, the results that we have obtained can be formulated as
follows
\begin{eqnarray}
\ss(L_n)\!\!&=&\!\left\{\begin{array}{cc}\!\!\!\!\!\!\!\!\!\!
\!\!\!\!\!\!\!\!\!\!\!\!u_1^nL_{n}
+nu_{1}^{n-1}v_1Y_{n}+w_{n}M_{n},\ \ \ n\neq0;\vs{6pt} \\
L_0+\frac{1}{2}(u_{1}w_{-1}-2u_1^{-2}v_{1}^2+u_{1}^{-1}w_{1})M_0,\ \
\ n=0;
\end{array}\right.\\
\ss(Y_n)\!&=&\!u_{1}^{n}e'_0Y_{n}+2nu_{1}^{n-1}v_1e'_0M_{n},\ \ \ n\in\Z;\\
\ss(M_n)\!\!\!&=&\!u_{1}^n{e'_0}^2M_{n},\ \ \ n\in\Z;\\
\ss
(C)\!\,&=&\!C+4(u_{1}w_{-1}+6u_1^{-2}v_{1}^2+u_{1}^{-1}w_{1}-u_1^{2}w_{-2}-u_1^{-2}w_{2})M_{0}.\label{eqN0013}
\end{eqnarray}
Replacing $\ss$ by $\psi_{{e'_0}^{-1}}\ss_{u_1^{-1}}\ss$ where
$\ss_{u_1^{-1}}$ and $\psi_{{e'_0}^{-1}}$ are defined in
(\ref{eqssa}) and (\ref{eqssa'}) and denoting
${e'_0}^{-1}u_1^{-1}v_1$ by $v'_1$, ${e'_0}^{-2}u_1^{-n}w_n$ by
$w'_n$, we can rewrite the above equations as follows
\begin{eqnarray}
\ss(L_n)\!\!&=&\left\{\begin{array}{cc}\!\!\!\!\!\!\!\!
\!\!\!\!\!\!\!\!\! \!\!\!L_{n}
+nv'_1Y_{n}+w'_{n}M_{n},\ \ \ n\neq0;\vs{6pt} \\
L_0+\frac{1}{2}(w'_{-1}-2{v'_{1}}^2+w'_{1})M_0,\ \ \ n=0;
\end{array}\right.\label{eqLateq11}\\
\ss(Y_n)\!&=&Y_{n}+2nv'_1M_{n},\ \ \ n\in\Z;\label{eqLateq21}\\
\ss(M_n)\!\!\!&=&M_{n},\ \ \ n\in\Z;\label{eqLateq31}\\
\ss
(C)\!\,&=&C+4(w'_{-1}+6{v'_{1}}^2+w'_{1}-w'_{-2}-w'_{2})M_{0}.\label{eqLateq41}
\end{eqnarray}
For $n\neq-1,0$, applying $\d$ to $[L_{1},L_{n}]=(n-1)L_{n+1}$, we
obtain
\begin{eqnarray*}
&\d\big([L_{1},L_{n}]\big)\!\!\!&=(n-1)\big(L_{n+1}
+(n+1)v'_1Y_{n+1}+w'_{n+1}M_{n+1}\big)\\
&&=(n-1)L_{n+1}+\big(n(n-\frac{1}{2})v'_1+(\frac{n}{2}-1)v'_1\big)Y_{n+1}\\
&&\ \ \ +\big(nw'_{n}+n(n-1){v'_1}^2-w'_{1}\big)M_{n+1}.
\end{eqnarray*}
Comparing the coefficients of $M_{n+1}$, one has
\begin{eqnarray}
&&(n-1)w'_{n+1}=nw'_{n}+n(n-1){v'_1}^2-w'_{1}.\label{eqwuv888}
\end{eqnarray}
In particular, one has
\begin{eqnarray}
w'_{-2}\!\!\!&=&\!\!\!w'_{-1}+3{v'_1}^2+\frac{w'_{-1}-w'_{1}}{2}=
w'_{1}+3{v'_1}^2+\frac{-3(w'_{1}-w'_{-1})}{2}.\label{eqwuv808}
\end{eqnarray}
For $n\neq-1,0$, applying $\d$ to $[L_{-1},L_{n+1}]=(n+2)L_{n}$, we
obtain
\begin{eqnarray*}
&\d\big([L_{-1},L_{n+1}]\big)\!\!\!&=(n+2)\big(L_{n}
+nv'_1Y_{n}+w'_{n}M_{n}\big)\\
&&=(n+2)L_{n}
+\big((n+1)(n+\frac{3}{2})v'_1-\frac{n+3}{2}v'_1\big)Y_{n}\\
&&\ \ \ +\big((n+1)w'_{n+1}-(n+1)(n+2){v'_1}^2+w'_{-1}\big)M_{n}.
\end{eqnarray*}
Comparing the coefficients of $M_{n}$, one has
\begin{eqnarray}
&&(n+1)w'_{n+1}=(n+2)w'_{n}+(n+1)(n+2){v'_1}^2-w'_{-1}.\label{eqwuv999}
\end{eqnarray}
In particular, one has
\begin{eqnarray}
w'_{2}=w'_{1}+3{v'_1}^2+\frac{w'_{1}-w'_{-1}}{2}.\label{eqwuv909}
\end{eqnarray}
For $n\neq-2,0,\pm1$, subtracting (\ref{eqwuv888}) from
(\ref{eqwuv999}), one has
\begin{eqnarray}
&&w'_{n+1}-w'_n=(2n+1){v'_1}^2+\frac{w'_{1}-w'_{-1}}{2}.\label{eqwuv888999}
\end{eqnarray}
According to the equations (\ref{eqwuv808}), (\ref{eqwuv909}) and
(\ref{eqwuv888999}), we obtain
\begin{eqnarray}
&&w'_{n}=w'_1+(n-1)(n+1){v'_1}^2+\frac{(n-1)(w'_{1}-w'_{-1})}{2},\ \
\ \forall\,\,n\in\Z^*.\label{eqwuv868999}
\end{eqnarray}
If we denote
$v'_1,\,\frac{w'_{1}-w'_{-1}}{2},\,w'_1-{v'_1}^2-\frac{w'_{1}-w'_{-1}}{2}$
respectively by $\a,\,\b,\,\g$, then for any $n\in\Z$, the equations
(\ref{eqLateq11})--(\ref{eqLateq41}) can be rewritten as
(\ref{eqLateq4200})--(\ref{eqLateq1200}). The lemma follows.\QED
\vskip7pt

 From (\ref{eqLateq4200})---(\ref{eqLateq1200}), one immediately see
\begin{eqnarray}
&&\delta_{\a,\b,\g}\delta_{\a',\b',\g'}=\delta_{\a+\a',\b+\b',\g+\g'+2\a\a'},\\
&&\delta_{\a,\b,\g}=\delta_{\a',\b',\g'}\ \Longleftrightarrow\
\a=\a',\,\,\b=\b',\,\,\g=\g'.
\end{eqnarray}
For any
\begin{equation}\label{c-infty}
b=(\cdots,b_{-2},b_{-1},b_1,b_2,\cdots),\
c=(\cdots,c_{-2},c_{-1},c_1,c_2,\cdots)\in\C^\infty,\end{equation}
we denote
\begin{equation}\label{xi---}\xi_{b,c}= {\rm
exp}\big(\mbox{$\sum\limits_{j\neq0}$}b_{j}{\rm
ad}Y_{j}+\!\mbox{$\sum\limits_{k\neq0}$}c_{k}{\rm
ad}M_{k}\big)\in\Inn(\LL).
\end{equation}
Then from the proofs of Lemmas \ref{Lemma 3.1}--\ref{Lemma 3.4}, we
see that every element $\si\in\Aut(\LL)$ can be uniquely expressed
as \big(the uniqueness can be seen from the fact that each data in
$(b,c,i,u,w,\a,\b,\g)$ is uniquely determined by $\si$\big)
\begin{equation}\label{unique-ex}
\si=\xi_{b,c}\eps^i\ss_u\psi_w\delta_{\a,\b,\g}\mbox{ for some
}(b,c,i,u,w,\a,\b,\g)\in\C^\infty\times\C^\infty\times\Z_2\times
\C^{*2}\times \C^3,
\end{equation}
where $\Z_2=\Z/2\Z$, and all notations can be found in
(\ref{eqssa})--(\ref{eqssa'}),
(\ref{eqLateq4200})---(\ref{eqLateq1200}) and (\ref{xi---}).
 Therefore there exists a
one to one correspondence
$f:\Aut(\LL)\to\C^\infty\times\C^\infty\times\Z_2\times
\C^{*2}\times \C^3$, \begin{equation}\label{corre}
f:\si\mapsto(b,c,i,u,w,\a,\b,\g).
\end{equation}
Suppose $\ss=\xi_{b,c}\eps^i\ss_{u}{\psi}_{{w}}\delta_{\a,\b,\g},\,
\ss'=\xi_{b',c'}\eps^{i'}\ss_{{u'}}\psi_{{w'}}\delta_{\a',\b',\g'}\in{\rm
Aut}(\LL)$. Let
$$\ss\ss'=\ss''=\xi_{b'',c''}\eps^{i''}\ss_{u''}\psi_{{w''}}\delta_{\a'',\b'',\g''}.$$
Then the data
$(b'',c'',i'',u'',w'',\a'',\b'',\g'')\in\C^\infty\times\C^\infty\times\Z_2\times
\C^{*2}\times \C^3$ can be determined by the following lemma.
\begin{lemm}\adddot
\label{last-lemma}
Under the above notations, the following relations hold:
\begin{eqnarray}
\!\!\!{w''}\!\!\!\!&=&\!\!\! {w}{w'},\label{eqnaq011}\\
\!\!\!i''\!\!\!&=&\!\!\! i+i',\label{eqnaq012}\\
\!\!\!{u''}\!\!\!&=&\!\!\! {u}^{(-1)^{i'}}{u'},\label{eqnaq013}\\
\!\!\!\g''\!\!\!&=&\!\!\! {w'}^{-2}\g+\g',\label{eqnaq014}\\
\!\!\!\a''\!\!\!&=&\!\!\! \frac{\a{w'}^{-1}+\a'}{2},\label{eqnaq015}\\
\!\!\!\b''\!\!\!&=&\!\!\! {w'}^{-2}\b+{\a'}^2+\b'+\g',\label{eqnaq016}\\
\!\!\!b''_{j}\!\!\!&=&\!\!\!b_{j} +(-1)^{i}{w}b'_{(-1)^{i}j}{u}^{(-1)^{i}j},\label{eqnaq017}\\
\!\!\!c''_{k}\!\!\!&=&\!\!\!c_{k}+{w}^2(-1)^{i}c'_{(-1)^{i}k}{u}^{(-1)^{i}k}
\!\!\!+2\a{w}^2kb'_{(-1)^{i}k}{u}^{(-1)^{i}k}\nonumber\\
&&\!\!\!\!\!\!\!-\mbox{$\sum\limits_{j\neq0}$}\!
\frac{(-1)^{i}{w}k^{-1}(k-j)(k-2j)\big({u}^{(-1)^{i}j}b'_{(-1)^{i}j}b_{k-j}
\!-\!{u}^{(-1)^{i}(k-j)}b_{j}b'_{(-1)^{i}(k-j)}\big)}{2}.\label{eqnaq018}
\end{eqnarray}
\end{lemm}
{\it Proof.} Using the equation $\ss\ss'(L_0)=\ss''(L_0)$ and
comparing the coefficients of $L_0$, $Y_j\,(j\neq0)$,
$M_k\,(k\neq0)$ and $M_0$, we can deduce (\ref{eqnaq012}),
(\ref{eqnaq017}), (\ref{eqnaq018}) and
\begin{eqnarray}
{w''}^2 \g''={w}^2\g+{w}^2{w'}^2\g'.\label{eqnaq020}
\end{eqnarray}
Applying both sides of $\ss\ss'=\ss''$ to $Y_0$ and comparing the
coefficients of $M_j\,(j\neq0)$, one can deduce (\ref{eqnaq011}).
Then by (\ref{eqnaq020}), we obtain (\ref{eqnaq014}).

Similarly, applying both sides of $\ss\ss'=\ss''$ to $M_1$ and
comparing the coefficients of $M_{(-1)^{i+i'}}$, one can deduce
(\ref{eqnaq013}).\

Using $\ss\ss'(Y_1)=\ss''(Y_1)$ and comparing the coefficients of
$M_{(-1)^{i+i'}}$, one has (\ref{eqnaq015}).

Finally, applying both sides of $\ss\ss'=\ss''$ to $L_1$ and
comparing the coefficients of $M_{(-1)^{i+i'}}$, one can deduce
(\ref{eqnaq016}).\QED \vskip6pt

Thus, we obtain the following theorem.
\begin{theo}\adddot
\label{last-theo}
Under the map $f$ defined in $(\ref{corre})$, the automorphism group
$\Aut(\LL)$ is ismorphic to the group
$\C^\infty\times\C^\infty\times\Z_2\times \C^{*2}\times \C^3$
$($where the elements in $\C^\infty$ is denoted as in
$(\ref{c-infty}))$, whose group multiplication is given by
$$(b,c,i,u,w,\a,\b,\g)\cdot(b',c',i',u',w',\a',\b',\g')=(b'',c'',i'',u'',w'',\a'',\b'',\g''),$$
where the data $(b'',c'',i'',u'',w'',\a'',\b'',\g'')$ is given by
Lemma $\ref{last-lemma}$.
\end{theo}
\begin{coro}\adddot
${\rm Inn}(\LL)={\rm Span}\big\{{\rm
exp}\big(aL_0+\mbox{$\sum_{i}b_iY_i$}
+\mbox{$\sum_{j\neq0}c_jM_j$}\big)\,|\,a,b_i,c_j\in\C\big\}$.
\end{coro}
{\it Proof.} The result follows immediately by noting from
(\ref{inn-group}) that $\Inn(\LL)$ is generated by ${\rm
exp\,}aL_0$, ${\rm exp\,}b_0Y_0$ and $\xi_{b,c}$ (which is defined
in (\ref{xi---})). \QED
\end{document}